\documentclass[12pt,reqno]{amsart}
\usepackage{enumerate, latexsym, amsmath, amsfonts, amssymb, amsthm, graphicx, color}
\def\pmod #1{\ ({\rm{mod}}\ #1)}
\def\Z{\Bbb Z}

\def\Q{\Bbb Q}

\def\bg{\bigg}
\def\({\bg(}
\def\){\bg)}

\def\sgn{{\rm sgn}}

\def\sgn{{\rm sgn}}
\def\Arg{{\rm Arg}}

\def\ov{\overline}

\def\ve{\varepsilon}

\theoremstyle{plain}
\newtheorem{theorem}{Theorem}

\newtheorem{lemma}{Lemma}

\newtheorem{conjecture}{Conjecture}
\theoremstyle{definition}
\newtheorem*{acknowledgment}{Acknowledgment}
\theoremstyle{remark}
\newtheorem{remark}{Remark}

\makeatletter
\@namedef{subjclassname@2010}{%
  \textup{2010} Mathematics Subject Classification}
\makeatother
 \vspace{4mm}

\begin{document}
 \baselineskip=17pt
\hbox{} {}
\medskip
\title[Applications of Lerch's theorem and permutations concerning quadratic residues]
{Applications of Lerch's theorem and permutations concerning quadratic residues}
\date{}
\author[Li-Yuan Wang and Hai-Liang Wu] {Li-Yuan Wang and Hai-Liang Wu}

\thanks{
\newline \indent Supported by the National Natural Science
Foundation of China (Grant No. 11571162).}

\address {(Li-Yuan Wang)  Department of Mathematics, Nanjing
University, Nanjing 210093, People's Republic of China}
\email{{\tt wly@smail.nju.edu.cn}}

\address {(Hai-Liang Wu) Department of Mathematics, Nanjing
University, Nanjing 210093, People's Republic of China}
\email{{\tt whl.math@smail.nju.edu.cn}}

\begin{abstract} Let $p$ be an odd prime. For each integer $a$ with $p\nmid a$, the famous Zolotarev's Lemma says
that the Legendre symbol $(\frac{a}{p})$ is
the sign of the permutation of $\Z/p\Z$ induced by multiplication by $a$. The extension of Zolotarev's result to
the case of odd integers was shown by Frobenius. After that, Lerch extended these to all positive integers. In this
paper we explore some applications of Lerch's result. For instance, we study permutations involving arbitrary $k$-th power residue modulo $p$ and primitive roots of a power of $p$. Finally, we discuss
some permutation problems concerning quadratic residues modulo $p$. In particular, we confirm some conjectures posed
by Sun.
\end{abstract}

\maketitle

\section{Introduction}
\setcounter{lemma}{0}
\setcounter{theorem}{0}
\setcounter{corollary}{0}
\setcounter{remark}{0}
\setcounter{equation}{0}
\setcounter{conjecture}{0}
For each integer $a$ and any positive integer $n$, throughout this paper, we let $\{a\}_n$ or $\ov {a}$ denote the
least nonnegative residue of $a$ modulo $n$ and $\sgn(\pi)$ denote the sign of a permutation $\pi$ over a finite set.
Following convention, $\sgn(\pi)=1$ if $\pi$ is even and $\sgn(\pi)=-1$ if $\pi$ is odd.

Let $p$ be an odd prime. For each integer $a$ with $p\nmid a$, the famous Zolotarev's Lemma \cite{Z} says
that the Legendre symbol $(\frac{a}{p})$ is
the sign of the permutation of $\Z/p\Z=\{\ov{0},\ov{1},...,\ov{p-1}\}$ induced by multiplication by $a$. This result
has many applications in modern number theory, and readers may consult \cite{BC,P} for more details.
Interestingly,  Zolotarev's Lemma can be generalized to all positive integers. Specifically, let $n$ be a positive
integer and $a$ is an integer prime with $n$. By elementary number theory we know that multiplication by $a$ induces
a permutation $\tau$ over $\Z/n\Z=\{\ov{0},\ov{1},...,\ov{n-1}\}$. Lerch obtained the following theorem which
determines the sign of $\tau$
(cf. \cite[Theorem 6.1]{BC}).

\begin{theorem}
Let notations as above, then  $$\sgn(\tau)=\begin{cases}(\frac{a}{n})&\mbox{if $n$ is odd},\\1&\mbox{if}\
n\equiv2\pmod4,\\(-1)^{\frac{a-1}{2}}&\mbox{if}\ n\equiv0\pmod4,\end{cases}$$
where $(\frac{\cdot}{\cdot})$ denotes Jacobi symbol.
\end{theorem}
It turns out that the case of even $n$ is very useful when we study permutation problems. In particular, we will use
this theorem and permutation polynomial to determine the sign of a permutation concerning arbitrary $k$-th power
residue
modulo an odd prime $p$.

Let $p$ be an odd prime and $k$ be a positive integer with $\gcd(p,k-1)=1$. Then it is easy to see that
$$\{1,2,3,\cdots,p-1\}=\{\{1^k\}_p,\{2^k\}_p,\{3^k\}_p,\cdots,\{(p-1)^k\}_p\}.$$
Noting also that $x^k\equiv 1\pmod p$ implies $x\equiv 1\pmod p$, we may therefore view
$$\{1^k\}_p,\{2^k\}_p,\{3^k\}_p,\cdots,\{(p-1)^k\}_p$$
as a permutation of $1,2,3,\cdots,p-1$ in this case. We denote this permutation by $\tau_{k,p}$. Here we mention that the
condition $\gcd(p,k-1)=1$ implies $x^k$ is a permutation polynomial over $\Z/p\Z=\{\ov{0},\ov{1},...,\ov{p-1}\}$. The
following theorem determines the sign of $\tau_{k,p}$.
\begin{theorem}\label{kth}
$$\sgn(\tau_{k,p})=\begin{cases}1&\mbox{if}\ p\equiv3\pmod4,\\(-1)^{\frac{k-1}{2}}&\mbox{if}\
p\equiv1\pmod4.\end{cases}$$
\end{theorem}

\begin{remark}
Let $p\equiv2\pmod3$ be an odd prime. Sun \cite{S} noticed that in this case $\sigma_3(\ov{k})=\ov{k^3}$ with $0\le
k\le p-1$
is a permutation on the set $\Z/p\Z=\{\ov{0},\ov{1},...,\ov{p-1}\}$. He conjectured that
$\sgn(\sigma_3)=(-1)^{(p+1)/2}$. This conjecture follows immediately from the above theorem.
\end{remark}

We now study permutations involving quadratic residues modulo an odd prime. Given an odd prime $p$,
let $1=a_1<a_2<...<a_{(p-1)/2}\le p-1$ be all quadratic residues modulo $p$ in the asscending order. It is easy to
see that  $a_1,a_2,...,a_{(p-1)/2}$ is a permutation over
$\{1^2\}_p, \{2^2\}_p,...,\{{(\frac{p-1}{2})}^2\}_p$. Let $\pi$ be this permutation. Sun discussed the sign of this
permutation. When
$p\equiv 3\pmod4$, he evaluated the product
$$\prod_{1\le j<k\le(p-1)/2}(\zeta_p^{j^2}-\zeta_p^{k^2})$$
by Galois theory,
where $\zeta_p=e^{2\pi i/p}$. With this result he managed to determine the sign of $\pi$ in case of $p\equiv 3 \pmod
4$:
$$\sgn(\pi)=\begin{cases}1&\mbox{if}\ p\equiv3\pmod8,\\
(-1)^{(h(-p)+1)/2}&\mbox{if}\ p\equiv7\pmod8.\end{cases}$$
In addition, he also studied some other permutations on quadratic
residues and posed some conjectures involving permutations of special forms. Readers
may consult \cite{S} for more details.

Inspired by Sun's work, we now consider the following sequences:
\begin{align}
A_0 &: a_1,a_2,...,a_{(p-1)/2},\\
A_1 &: \{1^2\}_p,\{2^2\}_p,...,\{(\frac{p-1}{2})^2\}_p,\\
A_2 &: \{2^2\}_p,\{4^2\}_p,...,\{(p-1)^2\}_p,\\
A_3 &: \{1^2\}_p,\{3^2\}_p,...,\{(p-2)^2\}_p,\\
A_4 &:  \{1\left(\frac{1}{p}\right)\}_p, \{2\left(\frac{2}{p}\right)\}_p,...,
\{\frac{p-1}{2}\left(\frac{(p-1)/2}{p}\right)\}_p.
\end{align}
where $\left(\frac{.}{p}\right)$ denotes Legendre symbol.
It is easy to see that $A_i \;(i=0,1,2,3)$ contains exactly all the quadratic residues modulo $p$. And $A_4$ does
only when $p\equiv
3\pmod 4$. If $A_i$ is a permutation of $A_j$, then we call this permutation $\sigma_{i,j}$. The following
theorem gives the sign of $\sigma_{2,1}$ and $\sigma_{3,1}$ for all odd primes.

\begin{theorem}\label{th213}
	$$\sgn(\sigma_{2,1})=\begin{cases}1&\mbox{if}\ p\equiv3\pmod4,\\(\frac{2}{p})&\mbox{if}\
p\equiv1\pmod4,\end{cases}$$
	and
	$$\sgn(\sigma_{3,1})=\begin{cases}-(\frac{2}{p})&\mbox{if}\ p\equiv3\pmod4,\\-1&\mbox{if}\
p\equiv1\pmod4.\end{cases}$$
\end{theorem}

\begin{remark}By the above theorem, it is easy to see that $\sgn(\sigma_{2,3})=(\frac{-2}{p})$.
\end{remark}

When $p\equiv 3\pmod4$, we determined the sign of $\sigma_{4,0}$ as follows:
\begin{theorem}\label{th40} Given an odd prime $p\equiv 3\pmod4$. Let $h(-p)$ denote the class number of
$\Q(\sqrt{-p})$,
	and let $\lfloor*\rfloor$ denote the floor function, then
	$$\sgn(\sigma_{4,0})=\begin{cases}(-1)^{\lfloor\frac{p+1}{8}\rfloor},&\mbox{if}\
	p\equiv3\pmod8,\\(-1)^{\lfloor\frac{p+1}{8}\rfloor+\frac{h(-p)+1}{2}}&\mbox{if}\ p\equiv7\pmod8.\end{cases}$$
\end{theorem}
\begin{remark}
	Here we note that combining Sun's result and the above theorem gives
	$$\sgn(\sigma_{4,1})=(-1)^{\lfloor\frac{p+1}{8}\rfloor}.$$
\end{remark}
In this line, Sun posed several conjectures, one of which states as follows. Letting $p$ be an odd prime and $k$ be
an integer, he defined $R(k,p)$ to be the unique $r\in\{0,1,..(p-1)/2\}$ with $k$ congruent
to $r$ or $-r$ modulo $p$ and
$$N_p:=\#\{(i,j): 1\le i<j\le \frac{p-1}{2}\ \text{and}\ R(i^2,p)>R(j^2,p)\},$$
where $\#S$ denotes the cardinality of a finite set $S$. He conjectured that
$N_p\equiv\lfloor\frac{p+1}{8}\rfloor\pmod2$ for every odd
prime $p$. Although we cannot prove this conjecture completely, we are able to obtain the following result.

\begin{theorem}\label{Np}
	Let notations be as above. For any prime $p\equiv3\pmod4$,
	$N_p\equiv\lfloor\frac{p+1}{8}\rfloor\pmod2$.
\end{theorem}

Let $p$ be an odd prime and let $A=\{1,2,...,(p-1)/2\}$. Sun \cite{S} defined a permutation $\tau_p$ as follows:
for each $k\in A$, $\tau_p(k)$ is the unique integer $k^*\in A$ with $kk^*\equiv\pm1\pmod p$. Sun \cite{S}
proved that
$\sgn(\tau_p)=-(\frac{2}{p})$. Here we give a simplier proof of this result by the aid of Lerch's theorem.

\begin{theorem}\label{sun}
Let notations be as above. Then $\sgn(\tau_p)=-(\frac{2}{p}).$
\end{theorem}

The proof of Theorem \ref{kth}--\ref{sun} will be given in the next section. In section 3, we turn to another kind of
permutation, which involves primitive roots.
\medskip
\maketitle
\section{Proof of Theorem \ref{kth}--\ref{sun}}
\setcounter{lemma}{0}
\setcounter{theorem}{0}
\setcounter{corollary}{0}
\setcounter{remark}{0}
\setcounter{equation}{0}
\setcounter{conjecture}{0}

\noindent{\it Proof of Theorem}\ 1.2. Let $g$ be a primitive root modulo $p$, then
$$\{1,2,3,\cdots,p-1\}=\{\{g^0\}_p,\{g^1\}_p,\{g^2\}_p,\cdots,\{g^{p-2}\}_p\}.$$ Since $$\tau_{k,p}(g^i)=g^{ki}\pmod
p$$ we see that  $\tau_{k,p}$ induces a permutation
$$\hat{\tau}(\ov{i})=\ov{ki}\pmod {p-1}$$ on the set $\Z/(p-1)\Z=\{\ov{0},\ov{1},...,\ov{p-2}\}$. It is easy to see
that $\hat{\tau}$ and $\tau_{k,p}$ have the same factorization. Hence
$$\sgn(\tau_{k,p})=\sgn(\hat{\tau})=\begin{cases}1&\mbox{if}\ p\equiv3\pmod4,\\(-1)^{\frac{k-1}{2}}&\mbox{if}\
p\equiv1\pmod4.\end{cases}$$
by Lerch's Theorem.
\qed

We now need the following lemma which originally appeared in \cite[pp.364--365]{SZ}.
\begin{lemma}\label{Lem B}
Let $p$ be a prime with $p\equiv 3\pmod4$. Then
$$\prod_{1\le i<j\le(p-1)/2}(i^2+j^2)\equiv (-1)^{\lfloor(p+1)/8\rfloor}\pmod p.$$
\end{lemma}

For convenience, we let $m=(p-1)/2$ throughout the remaining part of this section.

\noindent{\it Proof of Theorem}\ 1.3.  By definition, we have
\begin{eqnarray}
\sgn(\sigma_{2,1})
&\equiv&\prod_{1\le i<j\le \frac{p-1}{2}}^{}\frac{{(2j)}^{2}-{(2i)}^{2}}{j^2-i^2} \pmod p           \nonumber      \\
&=& \prod_{1\le i<j\le \frac{p-1}{2}}^{} 4 =4^{\frac{1}{2}\cdot \frac{p-1}{2}\cdot
\frac{p-3}{2}}\equiv\left(\frac{2}{p}\right)^{\frac{p-3}{2}} \pmod p           \nonumber      \\
&=&
\begin{cases} 1&\mbox{if}\ p\equiv3\pmod4\\ \left(\frac{2}{p}\right)&\mbox{if}\ p\equiv1\pmod4.\end{cases}
\end{eqnarray}
Now we calculate the sign of $\sigma_{3,1}$. By definition,

\begin{eqnarray}
\sgn(\sigma_{3,1})
&=&\prod_{1\le i<j\le m}^{}\frac{{(2j-1)}^{2}-{(2i-1)}^{2}}{j^2-i^2} \pmod p           \nonumber      \\
&\equiv& \prod_{1\le i<j\le m}^{} 4\cdot \frac{j+i-1}{j+i} \pmod p  \nonumber\\
&\equiv&2^{m\cdot  (m-1)}\cdot \frac{2}{m+1}\cdot \frac{4}{m+2}\cdot\cdots \frac{2m-2}{2m-1}\pmod p  \nonumber \\
&\equiv&2^{m\cdot (m-1)}\cdot 2^{m-1}\cdot \frac{(m-1)!\cdot m!}{ (2m-1)!}\pmod p    \nonumber\\
&\equiv& 2^{m^2-1}\cdot \frac{2}{p-1} \cdot(m!)^2  \pmod p \nonumber \\
&\equiv&-\left(\frac{2}{p}\right)^{\frac{p-1}{2}}\cdot (\frac{p-1}{2}!)^2  \pmod p .
\end{eqnarray}
This gives
\begin{equation*}
\sgn(\sigma_{3,1})=
\begin{cases} -1&\mbox{if}\ p\equiv1\pmod4\\ -\left(\frac{2}{p}\right)&\mbox{if}\ p\equiv3\pmod4\end{cases}
\end{equation*}\qed

\noindent{\it Proof of Theorem}\ 1.4. Here we assume $p\equiv 3\pmod 4$.
Since $a_1,a_2,...,a_m$ is the list of all the $(p-1)/2$ quadratic residues among  $1,..., p-1$ in the ascending
order, we only need to
count the number of ordered pairs $(i,j)$ with $1\leq i< j\leq m  $ and
$\{i\left(\frac{i}{p}\right)\}_p>\{j\left(\frac{j}{p}\right)\}_p$.
We denote this number by $s(p)$. Given any $1\leq i<m$, it is easy to check that if $\left(\frac{i}{p}\right)=1$ the
number of $j$ with
$1\leq i< j\leq m  $ and $\{i\left(\frac{i}{p}\right)\}_p>\{j\left(\frac{j}{p}\right)\}_p$ is zero and if
$\left(\frac{i}{p}\right)=-1$,
this number is $\frac{p-1}{2}-i$. Thus
\begin{align}
s(p)&=\sum\limits_{1\leq i\leq (p-1)/2}(\frac{p-1}{2}-i)\cdot \frac{1}{2}\left(1-\left(\frac{i}{p}\right)\right)
\nonumber\\
&=\frac{(p-1)(p-3)}{16}-\frac{p-1}{4}\sum\limits_{1\leq i\leq
(p-1)/2}\left(\frac{i}{p}\right)+\frac{1}{2}\sum\limits_{1\leq i\leq
(p-1)/2}i\left(\frac{i}{p}\right)
\end{align}
By Dirichlet's class number formula,
\begin{align}
-ph(-p)&=\sum\limits_{1\leq i\leq p-1}i\left(\frac{i}{p}\right) = \sum\limits_{1\leq i\leq (p-1)/2}
\left(i\left(\frac{i}{p}\right)+(p-i)\left(\frac{p-i}{p}\right) \right)  \nonumber\\
&=\sum\limits_{1\leq i\leq (p-1)/2} \left(2i\left(\frac{i}{p}\right)-p\left(\frac{i}{p}\right) \right).
\end{align}
This implies,
\begin{equation}
\sum\limits_{1\leq i\leq (p-1)/2} i\left(\frac{i}{p}\right)=\frac{1}{2}\left(-ph(-p)  +  p\sum\limits_{1\leq i\leq
(p-1)/2}\left(\frac{i}{p}\right)   \right)
\end{equation}
Thus
\begin{align}
s(p)&=\frac{(p-1)(p-3)}{16}-\frac{1}{4}ph(-p)+\frac{1}{4}\sum\limits_{1\leq i\leq (p-1)/2} \left(\frac{i}{p}\right)
\nonumber \\
&=\frac{(p-1)(p-3)}{16}-\frac{1}{4}ph(-p)+\frac{1}{4} \left(h(-p)-\left(\frac{2}{p}\right)\right)
\end{align}
The last equality follows from Dirichlet's class number formula in another form:
$$
h(-p)=\frac{1}{2-\left(\frac{2}{p}\right)}\sum\limits_{1\leq i\leq (p-1)/2} \left(\frac{i}{p}\right).
$$
When $p\equiv 3 \pmod 8$, letting $p=8k+3$ we get
\begin{align}
s(p)\equiv k  \pmod 2.
\end{align}
When $p\equiv 7 \pmod 8$, letting $p=8k+7$ we get
\begin{align}
s(p)\equiv k+1+\frac{h(-p)+1}{2} \pmod 2.
\end{align}
This gives
\begin{equation*}
s(p)\equiv
\begin{cases}\lfloor \frac{p+1}{8}\rfloor +\frac{h(-p)+1}{2} \pmod 2&\mbox{if}\ p\equiv3\pmod8,\\
 \lfloor \frac{p+1}{8}\rfloor\pmod 2&\mbox{if}\ p\equiv7\pmod8,\end{cases}
\end{equation*}
which completes the proof.\qed

\noindent{\it Proof of Theorem}\ 1.5.
Let $S$ be the set $\{1,2,...,\frac{p-1}{2}\}$ and $\tau :i\mapsto R(i^2,p) $ be a map from $S$ to itself. Since
$p\equiv 3\pmod 4$, it is
obvious that $\tau$ is a bijection thus also a permutation on $S$. Clearly,
$A_1=\{  \{1^2\}_p, \{2^2\}_p,...,\{{(\frac{p-1}{2})}^2\}_p    \}$ contains exactly all quadratic residues modulo
$p$. We define a map $f$
from $S$ to $A_1$ as follows: for each $k \in S$,$f(k)=\{ k^2\}_p$.
Thus
\begin{align}
\sgn(\sigma_{4,1})=\sgn(f\circ \sigma_{4,1}\circ f^{-1})=&\prod_{1\le i<j\le
\frac{p-1}{2}}^{}\frac{{j}^{4}-{i}^{4}}{j^2-i^2} \pmod p
\nonumber      \\
=&\prod_{1\le i<j\le \frac{p-1}{2}}^{}(j^2+i^2) \pmod p  \nonumber \\
= & (-1)^{\lfloor\frac{p+1}{8}\rfloor}  \pmod p.
\end{align}
The last equaltiy follows from Lemma \ref{Lem B}.\qed

\noindent{\it Proof of Theorem}\ 1.6.
For each integer $k$, recall that $\{k\}_p$ denotes the least nonnegative residue of $k$ modulo $p$. Then for each
$k\in A=\{1,2,...,(p-1)/2\}$, $\tau_p(k)=\{\ve_k k^{-1}\}_p$, where
$$\ve_k=\begin{cases}1&\mbox{if}\ 1\le\{k^{-1}\}_p\le (p-1)/2,\\-1&\mbox{otherwise}.\end{cases}$$
Let
$$B=\{\{1^2\}_p,\{2^2\}_p,...\{(\frac{p-1}{2})^2\}_p\}.$$
We define a map $f_1$ from $A$ to $B$ as follows:
for each $k\in A$, $f_1(k)=\{k^2\}_p$. Clearly, $f_1$ is a bijection. On the other hand,
let
\begin{align*}
A'&=\{\{\ve_11^{-1}\}_p,...,\{\ve_{\frac{p-1}{2}}(\frac{p-1}{2})^{-1}\}_p\},\\
B'&=\{\{1^{-2}\}_p,...,\{(\frac{p-1}{2})^{-2}\}_p\}.
\end{align*}
Then, we define a map $f_2$ from $A'$ to $B'$ as follows:
for each $\{\ve_kk^{-1}\}_p$,
$$f_2(\{\ve_kk^{-1}\}_p)=\{k^{-2}\}_p.$$
Clearly, $f_2$ is a bijection. Moreover, $f_2\circ \tau_p\circ f_1^{-1}$ is a permutation on $B$ with
$$f_2\circ \tau_p\circ f_1^{-1}(k^2)=\{k^{-2}\}_p.$$
It is easy to see that
$$\sgn(\tau_p)=\sgn(f_2\circ \tau_p\circ f_1^{-1}).$$
On the other hand, if we let $g$ be a primitive root of $p$, then
$$f_2\circ \tau_p\circ f_1^{-1}(g^{2l})=g^{-2l}.$$
Hence,
$f_2\circ \tau_p\circ f_1^{-1}$ induces a permutation $\pi_{-1}$ on
$\Z/(p-1)/2\Z=\{\ov{1},\ov{2},...,\ov{(p-1)/2}\}$.
For each $\ov{s}\in\Z/(p-1)/2\Z$, $\pi_{-1}(\ov{s})=\ov{-s}.$
Thus our theorem follows from Lerch's Theorem.

\maketitle
\section{Permutations Concerning Primitive Roots}
\setcounter{lemma}{0}
\setcounter{theorem}{0}
\setcounter{corollary}{0}
\setcounter{remark}{0}
\setcounter{equation}{0}
\setcounter{conjecture}{0}
In 2018, S. Kohl \cite{KLP} posed a permutation problem involving primitive roots of an odd prime on Mathoverflow.
Suppose that $p$ is an odd prime. Let $p$ be an odd prime, $\Z/p\Z=\{0,1,\ldots,p-1\}$ and $g$ is a primitive root of
$p$. Define
\begin{equation}
\sigma_g(b):=g^b
\end{equation}
for each $b\in\{1,\ldots,p-1\}$ and $\sigma_g(0)=0 $.
If we identify $\Z/p\Z$ with $\{0,1,\ldots,p-1\}$, we can view $\sigma_g$ as a permutation over $\Z/p\Z$. Let
$\mathcal{R}_{p}$ denote the set of all primitive roots of $p$. Kohl considered the sign of the permutation
$\sigma_g$ and raised the following
\begin{conjecture}Let notations be as above. Then
	
	{\rm (i)} If $p\equiv 1\pmod4$, then
	$$\#\{g\in\mathcal{R}_{p}: \sgn(\sigma_g)=1\}=\#\{g\in\mathcal{R}_{p}: \sgn(\sigma_g)=-1\}.$$
	
	{\rm (ii)} If $p\equiv 3\pmod4$, for each $g\in\mathcal{R}_{p}$ we have
	$$\sgn(\sigma_g)\equiv-\bigg(\frac{p-1}{2}\bigg)!\pmod{p}$$
	where $h(-p)$ denotes the class number of $\Q(\sqrt{-p})$.
\end{conjecture}
This conjecture was confirmed
by Ladisch and Petrov.
Throughout this section, we set $n=\phi(p^r)=p^{r-1}(p-1)$.
Now we discuss the sign of the permutation induced by the primitive roots of a power of an odd prime.
Given an odd prime $p$ and a positive integer $r$. Let $\mathcal{R}_{p^r}$ denote the set of all primitive roots of
$p^r$, and let
$1=b_1<b_2<...<b_n<p^r$ be the least nonnagetive reduced residue system modulo $p^r$ in ascending order, where
$n=\phi(p^r)=p^{r-1}(p-1)$.
Then for
each $g\in\mathcal{R}_{p^r}$, we define a permutation $\sigma_g$ on $\{b_1,..,b_n\}$ by
$$\sigma_g:\ b_i\mapsto g^i \pmod{p^r}.$$
We shall prove the following result.

\begin{theorem}\label{Th 1.5} Let notations be as above. Then
	
	{\rm (i)} If $p\equiv 1\pmod4$, then
	$$\#\{g\in\mathcal{R}_{p^r}: \sgn(\sigma_g)=1\}=\#\{g\in\mathcal{R}_{p^r}: \sgn(\sigma_g)=-1\}=n/2.$$
	
	{\rm (ii)} If $p\equiv 3\pmod4$, for each $g\in\mathcal{R}_{p^r}$, then
	$$\sgn(\sigma_g)=(-1)^{\frac{h-1}{2}},$$
	where $h(-p)$ denotes the class number of $\Q(\sqrt{-p})$.
\end{theorem}

\noindent{\it Proof of Theorem}\ 3.1 (i). When $p\equiv1\pmod4$, we have
$$\sigma_{g^{-1}}\circ\sigma_g^{-1}(g^i)=g^{-i}.$$
Thus $\sigma_{g^{-1}}\circ\sigma_g^{-1}$ induces a permutation $\pi_{-1}$ on
$\Z/n\Z=\{\ov{1}, \ov{2},...,\ov{n}\}$, where $\pi_{-1}(\ov{k})=\ov{-k}$. By Lerch's Theorem and
the fact that $p\equiv1\pmod4$, we obtain that $\sgn(\pi_{-1})=-1$. Thus
$\sgn(\sigma_{g^{-1}})\cdot \sgn(\sigma_{g})=-1$, which implies (i) of Theorem \ref{Th 1.5}.

(ii) Suppose $p\equiv 3\pmod4$. It follows from definition
$$\sgn(\sigma_g)=\prod_{1\le k<j\le n}\frac{\{g^j\}_{p^r}-\{g^k\}_{p^r}}{b_j-b_k}.$$
Thus we only need to determine
\begin{equation}\label{equation A}
\prod_{1\le k<j\le n}\frac{\{g^j\}_{p^r}-\{g^k\}_{p^r}}{b_j-b_k}\pmod p.
\end{equation}
We first consider the numerator. Let
\begin{equation}\label{equation B}
f(z)=\prod_{1\le k<j\le n} (z^j-z^k).
\end{equation}
Set $\zeta_n=e^{2\pi i/n}$, then
\begin{align*}
f(\zeta_n)^2=&(-1)^{\frac{n(n-1)}{2}}\prod_{1\le k\ne j\le n}(\zeta_n^j-\zeta_n^k)
=-1\cdot\prod_{1\le j\le n}\frac{z^n-1}{z-\zeta_n^j}\mid_{z=\zeta_n^j}\\
=&-1\cdot\prod_{1\le j\le n}n\zeta_n^{j(n-1)}=n^n.
\end{align*}
On the other hand, for each pair $(k,j)$ with $1\le k< j\le n$, it is easy to see that
$$\Arg(\zeta_n^j-\zeta_n^k)=\Arg(\zeta_n^{(j+k)/2}(\zeta_n^{(j-k)/2}-\zeta_n^{-(j-k)/2}))\equiv
\frac{j+k}{n}\pi+\frac{\pi}{2}\pmod
{2\pi},$$
where $\Arg(z)$ denotes the argument of complex number $z$. Thus
$$\Arg(f(\zeta_n))\equiv \sum_{1\le k<j\le n}(\frac{j+k}{n}\pi+\frac{\pi}{2})\equiv \frac{(3n+2)(n-1)}{4}\pi\pmod
{2\pi}.$$
Hence,
$$f(\zeta_n)=(-1)^{(3n+2)/4}n^{n/2}=p^{n(r-1)/2}(-1)^{(3n+2)/4}(p-1)^{n/2}.$$
Since $(\Z/p^r\Z)^{\times}$ is isomorphic to the group generated by $\zeta_n$, we have
\begin{equation}\label{numerator}
p^{-n(r-1)/2}\prod_{1\le k<j\le n}(\{g^j\}_{p^r}-\{g^k\}_{p^r}) \equiv (-1)^{(3n+2)/4+n/2}\pmod p.
\end{equation}
Now we consider the denominator. Since
$$\prod_{1\le k<j\le n}\frac{\{g^j\}_{p^r}-\{g^k\}_{p^r}}{b_j-b_k}=\pm1.$$
Thus we only need to determine
\begin{equation}\label{equation C}
p^{-n(r-1)/2}\prod_{1\le k<j\le n}(b_j-b_k)\pmod p.
\end{equation}
Note that
\begin{equation}\label{equation D}
p^{-n(r-1)/2}\prod_{1\le k<j\le n}(b_j-b_k)\equiv
\prod_{1\le i<j\le p-1}(j-i)^{p^{r-1}}\prod_{1\le i\ne j\le p-1}(j-i)^{\binom{p^{r-1}}{2}}\pmod p.
\end{equation}
It is known that
\begin{equation}\label{equation E}
\prod_{1\le i<j\le p-1}(j-i)\equiv (\frac{p-1}{2})!\cdot(-1)^{(p-3)/4}\pmod p.
\end{equation}
By \cite{M}, we know that
\begin{equation}\label{equation F}
(\frac{p-1}{2})!\equiv (-1)^{\frac{h(-p)+1}{2}}\pmod p,
\end{equation}
where $h(-p)$ denotes the class number of $\Q(\sqrt{-p})$.
Observe that
\begin{equation}\label{equation G}
\prod_{1\le i\ne j\le p-1}(j-i)=-1\times\prod_{1\le i<j\le p-1}(j-i)^2.
\end{equation}
Thus combining (\ref{equation D})-(\ref{equation G}), we obtain that
\begin{equation}\label{denominator}
p^{-n(r-1)/2}\prod_{1\le k<j\le n}(b_j-b_k)\equiv (-1)^{\frac{h(-p)+1}{2}+\frac{p-3}{4}+r+1}\pmod p.
\end{equation}
Our desired result follows from (\ref{numerator}) and (\ref{denominator}).
\qed

\acknowledgment\ We are grateful to Prof. Hao Pan for his helpful suggestions on the writting of this paper and
Prof. John Loxton for informing us of Lerch's theorem, which we proved as our theorem in the first version of this
paper. This research was supported by the National Natural Science Foundation of
China (Grant No. 11571162).

\end{document}